\documentclass[11pt]{article}
\usepackage{amsmath}
\usepackage{amssymb}
\usepackage{graphicx}

\newtheorem{theorem}{Theorem}[section]
\newtheorem{algorithm}[theorem]{Algorithm}

\newtheorem{proposition}[theorem]{Proposition}
\newtheorem{corollary}[theorem]{Corollary}
\newtheorem{lemma}[theorem]{Lemma}
\newtheorem{example}[theorem]{Example}

\def\qed{\hfill $\Box$\medskip}

\def\IC{{\bf C}}

\def\diag{{\rm diag}\,}
\def\tr{{\rm tr}\,}

\begin{document}
\openup .2\jot

\title{Product of two positive contractions}
\author{Chi-Kwong Li\thanks{Department of Mathematics, College of William and Mary,
Williamsburg,
VA 23187, USA. ckli@math.wm.edu}, Diane Christine Pelejo\thanks{Department of Mathematics,
College of William and Mary, Williamsburg, VA 23187, USA. dppelejo@wm.edu},
Kuo-Zhong Wang\thanks{Department of Applied Mathematics,
National Chiao Tung University, Hsinchu 30010, Taiwan.
kzwang@math.nctu.edu.tw}}
\date{}
\maketitle

\bigskip\centerline{\bf In memory of Professor Robert Thompson.}

\begin{abstract} Several characterizations are
given for a square matrix that can be written as the product of
two positive (semidefinite) contractions. Based on one of these characterizations,
and the theory of alternating projections, a Matlab program is written
to check the condition and construct the two positive contractions whose product equal
to the given matrix, if they exist.
\end{abstract}



\section{Introduction}

Let $M_n$ be the set of $n\times n$ complex matrices.
It is known that every matrix $A\in M_n$ with nonnegative determinant can be
written as the product of $k$ positive semidefinite matrices with $k \le 5$; see
\cite{B,CLS,W} and their references.
Moreover, characterizations are given
of matrices that can be written as the product of $k$
positive semidefinite matrices but not fewer for $k = 2, \dots, 5$.
In particular, a matrix $A$ is the product of two positive semidefinite matrices if
it is similar to a diagonal matrix with nonnegative diagonal entries.

In this paper, characterizations are given to
$A \in M_n$ which is a product of two
positive contractions, i.e., positive semidefinite matrices
with norm not larger than one. Evidently, if a matrix is the product of two
positive contractions, then it is a contraction
similar to a diagonal matrix with nonnegative diagonal entries.
However, the converse is not true. For example,
$A = \frac{1}{25} {\small \begin{pmatrix} 9 & 3  \cr 0 & 16 \cr\end{pmatrix}}$
is a contraction similar to $\diag(9,16)/25$ that is not a product of
two positive contractions as shown in \cite{LT}. In fact, the result in \cite{LT} implies that
if $A \in M_n$ is similar to a diagonal matrix with nonzero eigenvalues $a,b \in (0,1]$
then a necessary and sufficient condition for $A$ to be the product of two positive
contractions is:
$$\{\|A\|^2 - (a^2+b^2) + (ab/\|A\|)^2\}^{1/2} \le
|\sqrt{a} - \sqrt{b}| \sqrt{(1-a)(1-b)};$$
see Corollary \ref{cor}.
In particular, a matrix $A =
\begin{pmatrix}  a & p\cr 0 & b \cr\end{pmatrix}\in M_2$
is the product of two positive contractions if and only if
$a, b\in [0,1]$ and $|p| \le |\sqrt{a} - \sqrt{b}| \sqrt{(1-a)(1-b)}$.

In Section 2, we will present several characterizations of
a square matrix that can be written as the product of
two positive (semidefinite) contractions.
In Section 3, based on one of the characterizations in Section 2, we use
alternating projection method
to check the condition and construct the two positive contractions whose product equal
to the given matrix if they exist. Some numerical examples generated by Matlab are presented.

\section{Characterizations}

If $A$ is a product of two positive semidefinite contractions, then $A$ is similar to a diagonal
matrix with nonnegative eigenvalues with magnitudes bounded by
$\|A\| \le 1$. We will focus on such matrices in our characterization theorem.

It is known that a matrix $A$ is the product of two orthogonal projections if and only if it
is unitarily similar to a matrix which is the direct sum of $I_p\oplus 0_q$ and  matrices of the form
$$
\begin{pmatrix} a_j & \sqrt{a_j-a_j^2} \cr 0 & 0 \cr\end{pmatrix} \in M_{2}, \qquad
0<a_j<1\mbox{ for all }j= 1, \dots, m;
$$
see \cite{DLWWZ}.  Here we give another characterization which will be useful for our study.

\begin{proposition} \label{prop1}
Suppose $A$ is similar to $I_p \oplus 0_q\oplus \diag(a_1, \dots, a_m)$
with $a_1,\dots, a_m\in (0,1)$.
Then $A$ is the product of two orthogonal projections in $M_n$
if and only if $A$ is unitarily similar to $I_p \oplus A_1$ and there is
an $(n-p)\times m$ matrix $S$ of rank $m$ such that
 $A_1A_1^*S=A_1S=S\diag (a_1, \dots, a_m)$.
\end{proposition}

\it Proof. \rm For simplicity, we assume that $I_p$ is vacuous.
Suppose $A$ is the product of two orthogonal projections in $M_n$.
Let $D=\diag(a_1, \dots, a_m)$.
We may assume that $a_1\ge \cdots \ge a_m$.
There is a unitary $U$ such that $U^*AU=\left(
                                          \begin{array}{cc}
                                            D & \sqrt{D-D^2} \\
                                            0 & 0_m \\
                                          \end{array}
                                        \right)\oplus 0_{q-m}
$.
Let $U=[u_1\cdots u_n]$ and $U_m=[u_1\cdots u_m]$.
Hence, we have $AA^*S = AS = SD$ with $S=U_m$.

\medskip
Conversely, suppose $S$ satisfies
$AA^*S=AS=S\diag (a_1, \dots, a_m)$, and
has linearly independent columns $v_1, \dots, v_m$.
We may assume that $\|v_j\|=1$ for $1\le j\le m$ and
$\langle v_i,v_j\rangle=0$ if $a_i=a_j$ and $i\neq j$.
Since $AA^*$ is normal and $v_i$ is an eigenvector of $AA^*$ corresponding to
the eigenvalue $a_i$, $\langle v_i,v_j\rangle=0 $ for $a_i\neq a_j$.
Hence $S^*S=I_m$.
Now, we can find an orthonormal set $\{v_{m+1},\dots, v_{n}\}$ such that $V=[v_1, \dots, v_n]$ and $V^*AA^*V=D\oplus 0_{q}$.
Then $V^*AV$ is of the form $\left(
                            \begin{array}{cc}
                              D & B \\
                              0 & 0_{q} \\
                            \end{array}
                          \right)$,
where $B$ is an $m\times q$ matrix with $BB^*=D-D^2$.
From the QR factorization, $B$ can be written as $RQ$ with $Q$ unitary and $R$ lower triangular.
Let $V_1=I_m\oplus Q^*$.
Then $V_1^*V^*AVV_1=\left(
                      \begin{array}{cc}
                        D & R \\
                        0 & 0_{q} \\
                      \end{array}
                    \right)
$ and $RR^*=BQ^*QB^*=D-D^2$.
Hence $R=[\sqrt{D-D^2}\; 0_{m, (q-m)}]$, and we see that
$A$ is unitarily similar to the direct sum of $0_q$ and matrices of the form
$$
\begin{pmatrix} a_j & \sqrt{a_j-a_j^2} \cr 0 & 0 \cr\end{pmatrix} \in M_{2}, \qquad
j= 1, \dots, m.
$$
Hence $A$ is the product of two orthogonal projections.
\qed

Recall that $A \in M_n$ has a dilation $B \in M_N$ with $n < N$ if there is a unitary
$V \in M_N$ such that $A$ is the leading principal submatrix of $V^*BV$. For two Hermitian
matrices $X, Y \in M_n$, we write $X\ge Y$ if $X-Y$ is positive semidefinite.
In the next theorem, we present two characterizations for matrices which can be written as
the product of two positive contractions in terms of dilation and matrix inequalities.
We begin with the following observation.

\begin{lemma}
Suppose $A\in M_n$ is the product of two positive contractions.
Then $A$ is unitarily similar to a matrix of the form
$$
I_p\oplus
\begin{pmatrix}A_{11} & A_{12} \cr 0 & 0_{n-p-m} \cr\end{pmatrix},
$$
where $A_{11}\in M_m$ is similar to a diagonal matrix with the eigenvalues in $(0,1)$.
\end{lemma}

\it Proof. \rm
Obviously, the eigenvalues of $A$ are in $[0,1]$.
From [2, Proposition $3.1(d)$], we have
$$
A\cong \left(
         \begin{array}{ccc}
           I_p & B_1 & B_2 \\
           0 & A_{11} & A_{12} \\
           0 & 0 & 0_{n-p-m} \\
         \end{array}
       \right),
$$
 where $A_{11}\in M_m$ is an upper block triangular matrix such that the diagonal blocks are scalar matrices corresponding to distinct scalars, $1>\lambda_1>\cdots>\lambda_k>0$.
 Since $\|A\|\le 1$, $B_1$ and $B_2$ are zero matrices.
 By [2, Proposition $3.1(c)$ and $(d)$], $A_{11}$ is similar to a diagonal matrix, and the desired conclusion follows.
 \qed

\begin{theorem} \label{thm1}
Suppose $A = I_p \oplus
\begin{pmatrix}A_{11} & A_{12} \cr 0 & 0_{n-p-m} \cr\end{pmatrix} \in M_{n}$
such that $A_{11}\in M_m$ is similar to
$D\equiv\diag(a_1, \dots, a_m)$ with $1 > a_1 \ge \cdots \ge a_m > 0$.
The following conditions are  equivalent.

\begin{itemize}
\item[{\rm (a)}] $A$ is the product of two positive contractions.
\item[{\rm (b)}] $A$ has a dilation $\tilde T \in M_{n+2m}$, which is the product of
two orthogonal projections and has the same rank and eigenvalues of $A$.
Equivalently, there are matrices $R, C \in M_m$ such that
$$\tilde T = I_p \oplus  {\small \begin{pmatrix} A_{11} & A_{12} & 0 & A_{11}C \cr
 0 & 0_{n-p-m} & 0 & 0 \cr
RA_{11} & RA_{12} & 0_m & RA_{11}C
\cr 0 & 0 & 0 & 0_m
\cr\end{pmatrix}}
\in M_{n+2m}$$
is the product of two orthogonal projections.
\item[{\rm (c)}]
There is an invertible contraction $U_{11}\in M_m$ satisfying
$$A_{11}U_{11}=U_{11}D\quad \mbox{ and }\quad U_{11}DU_{11}^* \ge A_{11}A_{11}^* + A_{12}A_{12}^*.$$
\end{itemize}
Moreover, if condition {\rm (c)} holds,
we have $A = (I_p \oplus P)(I_p \oplus Q)$ for the positive contractions
$$P = \begin{pmatrix} U_{11}U_{11}^*& 0 \cr 0 & 0_{n-p-m} \cr\end{pmatrix}
\quad \hbox{ and } \quad
Q = \begin{pmatrix} (U_{11}^*)^{-1}DU_{11}^{-1}& (U_{11}U_{11}^*)^{-1}A_{12} \cr
A_{12}^*(U_{11}U_{11}^*)^{-1} &  A_{12}^*(U_{11}DU_{11}^*)^{-1}A_{12}\cr\end{pmatrix}.$$
\end{theorem}

\it Proof. \rm
For simplicity, we can assume that $I_p$ is vacuous because the matrix
$A$ is the product of two positive
contractions if and only if each of the two positive
contractions is a direct sum of $I_p$ and a positive contraction in $M_{n-p}$.

First we establish the equivalence of (a) and (b).
If (a) holds, then
$A = P Q$, where $P, Q$ are two positive contractions. Then
$${\footnotesize
\tilde P = \begin{pmatrix} P & \sqrt{P-P^2} & 0 \cr  \sqrt{P-P^2} & I_n - P& 0\cr
0&0&0_n\cr\end{pmatrix}} \quad \hbox{ and } \quad
{\footnotesize
\tilde Q = \begin{pmatrix} Q &  0 &  \sqrt{Q-Q^2} \cr  0 & 0_n & 0\cr
\sqrt{Q-Q^2} & 0 & I_n - Q\cr\end{pmatrix}}$$
are orthogonal projections such that
$$\tilde P\tilde Q =
{\small \begin{pmatrix} PQ  & 0 & P\sqrt{Q-Q^2}\cr
\sqrt{P-P^2}Q & 0_n& \sqrt{(P - P^2)(Q-Q^2)} \cr
0&0&0_n\cr
\end{pmatrix}}.$$
Let $Y =  \sqrt{Q^+-Q^+Q}$ and $X = \sqrt{P^+ - P^+P}$,
where $P^+,Q^+$ is the Moore-Penrose inverses of $P$ and $Q$.
(Recall that for a Hermitian matrix
$H = \sum_{j=1}^\ell \lambda_j \xi_j \xi_j^*\in M_n$
with nonzero  eigenvalues $\lambda_1, \dots, \lambda_\ell$ and orthonormal
eigenvectors $\xi_1, \dots, \xi_\ell$, its Moore-Penrose inverse $H^+$ is
$\sum_{j=1}^\ell \lambda_j^{-1} \xi_j \xi_j^*$.)
Let
$$T = {\small \begin{pmatrix}A & 0 & AY \cr X^*A & 0_n & X^*AY \cr 0 & 0 & 0_n \cr\end{pmatrix}.}$$
The rows of the matrix $X^*A$ lie in the row space of $[A_{11} A_{12}]$
and the columns of $AY$ lie in the column space of $A_{11}$. So, there is unitary matrix
of the form  $U = I_n \oplus U_1 \oplus U_2$ with $U_1, U_2 \in M_n$
such that
$$U^*TU = {\small \begin{pmatrix}
A_{11} & A_{12} & 0_m & 0_{m,n-m} & A_{11}C & 0_{m,n-m}\cr
0_{n-m,m} & 0_{n-m} & 0_{n-m,m} & 0_{n-m} & 0_{n-m,m} & 0_{n-m}\cr
RA_{11} & RA_{12} & 0_{m} & 0_{m,n-m} & RA_{11}C& 0_{m,n-m}\cr
0_{n-m,m} & 0_{n-m} & 0_{n-m,m} & 0_{n-m} & 0_{n-m,m} & 0_{n-m}\cr
0_{n,m} & 0_{n,n-m} & 0_{n,m} & 0_{n,n-m} & 0_{n,m} & 0_{n,n-m}\cr
\end{pmatrix}}.$$
Thus,
$$\tilde T =
{\small  \begin{pmatrix}
A_{11} & A_{12} & 0 & A_{11}C\cr
0 & 0_{n-m} & 0 & 0 \cr
RA_{11} & RA_{12} & 0_{m} & RA_{11}C \cr
0 & 0 & 0 & 0_{m}\cr
\end{pmatrix} \in M_{n+2m}}$$
has the same rank and eigenvalues as the leading submatrix $A$.
Thus, condition (b) holds.

Conversely, suppose (b) holds.
and $\tilde T$ is the product of two orthogonal projections
$\tilde P = VV^*$ and $\tilde Q = WW^*$ with $V \in M_{n+2m,r}, W \in M_{n+2m,s}$
such that $V^*V = I_r$ and $W^*W = I_s$.
Evidently, $\tilde T$ has rank $m$. So,
$$V^*W = Y\begin{pmatrix}K & 0 \cr 0 & 0_{(r-m),(s-m)} \cr\end{pmatrix}Z^*$$
such that $Y \in M_r$, $Z \in M_s$ are unitary and $K \in M_m$ is a diagonal matrix with positive diagonal entries.
Let $Y = [Y_1|Y_2], Z = [Z_1|Z_2] $ be such that
$Y_1 \in M_{r,m}, Z_1 \in M_{s,m}$.
Note that
$$Y_1^*V^*WZ_1=Y_1^*[Y_1|Y_2]\begin{pmatrix}K & 0 \cr 0 & 0_{(r-m),(s-m)} \cr\end{pmatrix}[Z_1|Z_2]^*Z_1=K.
$$
Furthermore,
$$\tilde V = VY_1 =  \begin{pmatrix}V_1 \cr V_2 \cr V_3\cr\end{pmatrix} \quad  and \quad
\tilde W =WZ_1= \begin{pmatrix}W_1 \cr W_2 \cr W_3\cr\end{pmatrix},$$
where $V_1, W_1$ are $n\times m$, $V_2,V_3, W_2, W_3 \in M_m$.
Then
$$\tilde V\tilde V^* \tilde W \tilde W^*
= VY_1Y_1^*V^* WZ_1Z_1^*W^* = VY_1KZ_1^*W^* = VV^*WW^* = \tilde T.$$
Now, the last $m$ rows of $\tilde T$ and the $(n+1)$st, $\dots, (n+m)$th columns of $\tilde T$
are zero. Thus,
$$V_3 \tilde V^* \tilde W \tilde W^* =
V_3  K \tilde W^* = 0_{m, (n+2m)} \quad \hbox{ and } \quad
\tilde V\tilde V^* \tilde W W_2^* = \tilde V  K W_2^* = 0_{(n+2m), m}.$$
Because $ K\tilde W^*$ has full row rank and
$\tilde V K$ has full column rank, we see that
$V_3 = 0_m$ and $W_2 = 0_m$. Consequently, $A = V_1V_1^*W_1W_1^*$ is the product of
two positive contractions $V_1V_1^*$ and $W_1W_1^*$.

Next, we prove the equivalence of conditions (b) and (c).
Suppose (b) holds, and
$$\tilde T =
 \begin{pmatrix}
A_{11} & A_{12} & 0 & A_{11}C\cr
0 & 0_{n-m} & 0 & 0 \cr
RA_{11} & RA_{12} & 0_{m} & RA_{11}C \cr
0 & 0 & 0 & 0_{m}\cr
\end{pmatrix} \in M_{n+2m}$$
has the same rank and eigenvalues as the leading submatrix $A$.

Now, assume that $U = (U_{ij})_{1 \le i \le 4, 1\le j \le 3} \in M_{n+2m}$ is unitary
with $U_{11}, U_{12} \in M_m, U_{13} \in M_{m, n}$
and $U_{31}, U_{41} \in M_m, U_{21} \in M_{n-m,m}$
such that
$$U^*\tilde T U = \begin{pmatrix} D &  \sqrt{D-D^2} & 0\cr
 0 & 0_m & 0\cr
 0 & 0 & 0_{n}\cr \end{pmatrix}.$$
 Now,
$$\begin{pmatrix}
A_{11} U_{11} + A_{12} U_{21} \cr 0_{n-m,m} \cr RA_{11}U_{11} + RA_{12}U_{21} \cr 0_m\cr
\end{pmatrix}
= \tilde T \begin{pmatrix} U_{11} \cr U_{21} \cr U_{31} \cr U_{41} \cr \end{pmatrix} =
 \begin{pmatrix} U_{11} \cr U_{21} \cr U_{31} \cr U_{41} \cr \end{pmatrix}D.$$
 It follows that $U_{21}, U_{41}$ are zero matrices.
Furthermore,
$$A_{11} U_{11} = U_{11}D, \quad RA_{11}U_{11} = U_{31}D.$$
Thus, $RU_{11}D = U_{31}D$ so that $RU_{11} = U_{31}$.
If $x \in \IC^m$ satisfies
$U_{11} x = 0$, then
$$x = (U_{11}^* \ U_{31}^*)
\begin{pmatrix} U_{11} \cr U_{31} \cr \end{pmatrix} x = U_{11}^*(I_m+R^*R)U_{11}x  = 0.$$
Hence, $U_{11} \in M_m$ has linearly independent columns, i.e., $U_{11}$ is invertible.

Next, observe that
$$\tilde T\tilde T^*U = U \begin{pmatrix} D &  0 & 0\cr
 0 & 0_m & 0\cr
 0 & 0 & 0_{n}\cr \end{pmatrix}.$$
 So,
 $$(A_{11}A_{11}^* + A_{12}A_{12}^* + A_{11}CC^*A_{11}^*)(I_m+R^*R)U_{11} = U_{11}D,$$
and hence
\begin{equation}\label{CC}
(A_{11}A_{11}^* + A_{12}A_{12}^* + A_{11}CC^*A_{11}^*) = U_{11}DU_{11}^*,
\end{equation}
because
\begin{equation}\label{RR}
I_m = U_{11}^*U_{11} + U_{31}^*U_{31} = U_{11}^*(I_m+R^*R)U_{11} = (I_m+R^*R)U_{11}U_{11}^*.
\end{equation}
So, $R$ and $C$ exist if and only if there is  a contraction $U_{11}\in M_m$
satisfying
$$A_{11}U_{11} = U_{11}D \quad \hbox{ and } \quad
U_{11}DU_{11}^* \ge A_{11}A_{11}^* + A_{12}A_{12}^*.$$

Conversely, suppose $(c)$ holds.
Then there exist $R$ and $C$ satisfying (\ref{CC}) and (\ref{RR}). Let
$$\tilde U = {\small \begin{pmatrix}U_{11} \cr 0_{n-m,m} \cr RU_{11} \cr 0_m\cr\end{pmatrix}}.$$
Then $\tilde U$ has rank $m$ and the matrix
$\tilde T$ in condition (b) satisfies $\tilde T \tilde T^*\tilde U=\tilde T \tilde U = \tilde UD$.
By Proposition \ref{prop1}, we see that $\tilde T$ is the product of two orthogonal projections.

\medskip
To verify the last statement, note that $A_{11}U_{11} = U_{11}D$
so that $A_{11} = U_{11}DU_{11}^{-1}$. Hence,
$$PQ = {\small \begin{pmatrix}U_{11}DU_{11}^{-1} &
A_{12} \cr 0 & 0_{n-m} \cr\end{pmatrix}}
= {\small \begin{pmatrix}A_{11} &
A_{12} \cr 0 & 0_{n-m} \cr\end{pmatrix}},$$
and
$Q = ZZ^*$ with $Z = \begin{pmatrix}
(U^*_{11})^{-1}D^{1/2} \cr A_{12}^*(U_{11}^*)^{-1}D^{-1/2}\cr\end{pmatrix}$
so that
\begin{eqnarray*}Z^*Z &=&
D^{1/2}U_{11}^{-1}(U_{11}^*)^{-1}D^{1/2} +
D^{-1/2}U_{11}^{-1}A_{12}A_{12}^*(U_{11}^*)^{-1}D^{-1/2}  \\
&=& D^{-1/2}U_{11}^{-1} (A_{11}A_{11}^*+A_{12}A_{12}^*) (U_{11}^*)^{-1}D^{-1/2} \\
&\le& D^{-1/2}U_{11}^{-1} (U_{11}DU_{11}^*) (U_{11}^*)^{-1}D^{-1/2} = I_m.
\end{eqnarray*}
This shows that $Z$ is a contraction and hence so is $Q$.
\qed

\medskip
As pointed out by the referee,
from Theorem \ref{thm1} one can deduce the following corollary,
which can be viewed as a 2-variable
generalization of the fact that every positive contraction can be dilated 
to an orthogonal projection; see \cite[Problem 222(b)]{Halmos}.

\begin{corollary}
If $A\in M_n$ is the product of two positive contractions, 
then $A$ can be dilated to a product of two projections on $\mathbb{C}^{n+2m}$, 
where $m$ equals the number of eigenvalues of $A$ which are not equal to $0$ or $1$.
\end{corollary}

It is not easy to check the existence of the matrices $R, C \in M_m$ in condition (b),
and the existence of $U_{11}$ in condition (c) of Theorem \ref{thm1}.
We refine condition (c) to get Theorem \ref{thm2} below so that one can use
computational techniques such as positive semidefinite programming or alternating projection
methods to check the condition.
In Section 3, we will develop Matlab programs
using an alternating projection method based on
Theorem \ref{thm2}  to check whether a matrix can be written as
the product of two positive semidefinite contractions, and construct
them if they exist.

\begin{theorem} \label{thm2}
Let $A \in M_n$ be unitarily similar to $I_p \oplus 0_q \oplus
{\small \begin{pmatrix}A_{11} & A_{12} \cr 0 & 0_{n-p-q-m} \cr\end{pmatrix}}$,
where $A_{11}\in M_m$ such that $A_{11}$ is diagonalizable with distinct
eigenvalues $\alpha_1 > \cdots > \alpha_k$ in $(0,1)$ with multiplicities
$m_1, \dots, m_k$, respectively.
Suppose $V= [V_1 \ \cdots \ V_k] \in M_m$ is an invertible matrix
such that the columns of the $n\times m_j$ matrix
$V_j$ form an orthonormal basis for the null space of $A_{11}-\alpha_j I_m$,
for $j = 1, \dots, k$, i.e., $A_{11}V = VD$, where
$D = \alpha_1 I_{m_1} \oplus \cdots \oplus \alpha_k I_{m_k}$
and $V_j^*V_j = I_{m_j}$ for $j=1,\ldots,k$.
Then $A$ is the product of two positive contractions if and only if
there is a block diagonal matrix
$\Gamma = \Gamma_1 \oplus \cdots \oplus \Gamma_k\in
M_{m_1} \oplus \cdots \oplus M_{m_k}$ satisfying
\begin{equation}\label{eq}
D^{1/2}V^*(A_{11}A_{11}^* + A_{12}A_{12}^*)^{-1}VD^{1/2}\ge \Gamma \ge V^*V.
\end{equation}
\end{theorem}

\it Proof. \rm Suppose $A_{11}V = VD$ as asserted. Then  $U$ satisfies
$A_{11}U = UD$ if and only if $U = VL$ for some block matrix
$L = L_1 \oplus \cdots \oplus L_k \in M_{m_1} \oplus \cdots \oplus M_{m_k}$.
One readily checks that condition (c) in Theorem  \ref{thm1} reduces to the
existence of $\Gamma = (LL^*)^{-1}$.
\qed

By Theorem \ref{thm2}, we can deduce the following corollary.
The first part of the corollary was obtained in \cite[Lemma 2.1]{LT}
by some rather involved arguments. The second part of the corollary is a proof
of a comment in our introduction.

\begin{corollary} \label{cor} Let
$A = \begin{pmatrix} a & p \cr 0 & b \cr\end{pmatrix}$
 with $a, b \in [0,1]$. Then $A$ is the product of two positive
contractions if and only if
\begin{equation}\label{ineq}
|p| \le |\sqrt{a}-\sqrt{b}| \sqrt{(1-a)(1-b)}.\end{equation}
Consequently,
if $B \in M_n$ is similar to a diagonal matrix with nonzero eigenvalues $a,b \in (0,1]$
then a necessary and sufficient condition for $A$ to be the product of two positive
contractions is:
$$\{\|B\|^2 - (a^2+b^2) + (ab/\|B\|)^2\}^{1/2} \le
|\sqrt{a} - \sqrt{b}| \sqrt{(1-a)(1-b)}.$$
\end{corollary}

\it Proof. \rm   Case 1. $a = b$. If $A$ is the product of two positive contractions,
then $A$ is similar to a diagonal matrix so that $p = 0$, and inequality (\ref{ineq}) holds.
If inequality (\ref{ineq}) holds, then $p = 0$, and $A = aI_2$ is the product of positive
contractions $I_2$ and $aI_2$.

Case 2. $a \ne b$.
We focus on the non-trivial case that $a,b \in (0,1)$, $a\neq b$ and $p \ne 0$.
One sees that $V$ in Theorem \ref{thm2} can be chosen to be
$\begin{pmatrix}1 & p/\gamma \cr 0 & (b-a)/\gamma \cr\end{pmatrix}$
with $\gamma = \sqrt{(a-b)^2+p^2}$ so that up to
diagonal congruence  we have
$$V^*V = \begin{pmatrix}1 & p/\gamma \cr p/\gamma & 1 \cr\end{pmatrix}.$$
We need to find a diagonal matrix $\Gamma = \diag(d_1,d_2)$ with $d_1, d_2 \ge 0$ such that
$\Gamma - V^*V \ge 0$ and $VV^* - \diag(ad_1, bd_2) \ge 0$.
Thus, we want
$$(d_1-1)(d_2-1) \ge p^2/\gamma^2, \qquad (1-d_1a)(1-d_2b) \ge p^2/\gamma^2.$$
We consider the maximum values for
$$f(d_1,d_2) = (d_1-1)(d_2-1)$$
subject to the condition of
$$g(d_1,d_2) = (d_1-1)(d_2-1) - (1-d_1a)(1-d_2b) = 0.$$
Consider the Lagrangian function
$L(d_1,d_2,\mu) = f(d_1,d_2) - \mu g(d_1,d_2)$.
$$0 = L_{d_1}(d_1,d_2,\mu) = (d_2-1) - \mu[(d_2-1) + a (1-d_2b)]$$
and
$$0 = L_{d_2}(d_1,d_2,\mu) = (d_1-1) - \mu[(d_1-1) + b (1-d_1a)].$$
Thus,
$$(1-\mu)^2(d_1-1)(d_2-1) = \mu^2 ab(1-d_1a)(1-d_2b).$$
Because $(d_1-1)(d_2-1) = (1-d_1a)(1-d_2b)$, we see that
$(1-\mu)^2 = \mu^2 ab$, and thus, $\mu = (1+\sqrt{ab})^{-1}$.
Here, we use the root satisfying $1-\mu > 0$.
Solving $d_1$ and $d_2$, we get
$$(d_1-1)(d_2-1) = (1-a)(1-b)/(1+\sqrt{ab})^2.$$
Furthermore, $(d_1-1)(d_2-1) \ge p^2/\gamma^2$
if and only if
$$p^2 \le (a-b)^2(1-a)(1-b)/(\sqrt{a}+\sqrt{b})^2 = (\sqrt{a} - \sqrt{b})^2(1-a)(1-b).$$

For the last assertion, note that if $B$ satisfies the given assumption, then
$(B-aI)(B-bI) = 0$, and $B$ is unitarily similar to the direct sum of
$aI_p \oplus bI_l$ and matrices of the form
$B_j = \begin{pmatrix} a & p_j \cr 0 & b \cr\end{pmatrix}$,
where $p_1 \ge \cdots \ge p_k > 0$, for $j = 1, \dots, k$.
By Theorem 1.1 in \cite{LT}, $B$ is a product of two positive contractions if and only if
$$\|\diag(p_1, \dots, p_k)\| =
|p_1| \le  |\sqrt{a}-\sqrt{b}| \sqrt{(1-a)(1-b)}.$$
It is easy to check that $\|B\| = \|B_1\|$ and
$$\|B_1\|^2+ (ab/\|B_1\|)^2 - (a^2+b^2) = \tr(B_1^*B_1) - (a^2+b^2) = p_1^2.$$
The assertion follows. \qed

\section{Alternating projections and numerical examples}

In Theorem \ref{thm2}, if $A_{11}$ has distinct eigenvalues, then one only
needs to search for a diagonal matrix satisfying the condition.
However, there is no guarantee that there is a diagonal matrix $\Gamma$ satisfying
the condition in general as shown in the following example.

\begin{example} \rm
Let $D = \diag(0.15,0.15, 0.2)$,
$A = \begin{pmatrix} A_{11} & A_{12} \cr 0_3 & 0_3\cr\end{pmatrix}$
with
$$A_{11} ={\small \begin{pmatrix}
 0.1500 &        0 &        0\cr
      0 &    0.1500&   0.0375\cr
      0 &        0 &  0.2000\cr
\end{pmatrix}},$$
and
$$
A_{12} = \{UDU^*-A_{11}A_{11}^*\}^{1/2} =
{\small \begin{pmatrix}
   0.3571&         0&          0\cr
        0&    0.3215&    0.1070\cr
        0&    0.1070&    0.1689\cr
\end{pmatrix}},$$
where
$$U = V R = {\small
\begin{pmatrix}
1&         0&         0\cr
0&    5/\sqrt{40}&  3/\sqrt{40} \cr
0&         0&    4/\sqrt{40}\cr
\end{pmatrix}},$$
with
$$
V =  {\small
\begin{pmatrix}
1/\sqrt{2}&  1/\sqrt{2} &  0\cr
1/\sqrt{2}&  -1/\sqrt{2}&  3/5\cr
0&         0&    4/5\cr
\end{pmatrix}} \quad \hbox{and} \quad
R = {\small
\begin{pmatrix}
1/\sqrt{2}&  1/\sqrt{2} &  0\cr
1/\sqrt{2}&  -1/\sqrt{2}&  0\cr
0&         0&    1\cr
\end{pmatrix}}
{\small
\begin{pmatrix}
1&  0 &  0\cr
0 & 5/\sqrt{40}&  0\cr
0&         0&    5/\sqrt{40}\cr
\end{pmatrix}}.$$
\end{example}
\noindent
Then $A_{11}V = VD$, $A_{11}U = UD$,  and $U$ is a contraction such that
$UDU^* = A_{11}A_{11}^*+A_{12}A_{12}^*$.
There is no $\Gamma = \diag(\mu_1, \mu_2, \mu_3)$ such that
$$M = D^{1/2}V^*(A_{11}A_{11}^* + A_{12}A_{12}^*)^{-1}VD^{1/2}
= {\small\begin{pmatrix}
  1.3&   -0.3&        0\cr
   -0.3&    1.3&      0\cr
      0&      0&    1.6\cr
\end{pmatrix}}
\ge \Gamma
$$
and
$$
\Gamma \ge V^*V = {\small
\begin{pmatrix}
 1.0000&         0&    0.4243\cr
         0&    1.0000&   -0.4243\cr
    0.4243&   -0.4243&    1.000\cr
\end{pmatrix}}$$
because $\mu_1, \mu_2 \in (1, 1.3)$ so that the leading $2\times 2$ principal submatrix
$M - \Gamma$ cannot be positive semidefinite. Hence, $A$ is not the product
of two positive contractions.\qed

By Theorem \ref{thm2}, one can use positive semidefinite (PSD) programming  to check whether
there exists $\Gamma$ satisfying (\ref{eq}). However, standard PSD programming uses
dual program to check the feasibility, and does not seem to be effective in checking
the result. For example, we use the the SDP mode of cvx program from \texttt{http://cvxr.com/cvx/}, and
it fails to detect the result even for $A \in M_2$.

We turn to alternating projection method; for example see \cite{EscRayd}.
Suppose $A \in M_n$ is a contraction matrix unitarily similar to $I_p \oplus 0_q \oplus
\begin{pmatrix}A_{11} & A_{12} \cr 0 & 0_{n-p-q-m} \cr\end{pmatrix}$
and $V \in M_m$ is an invertible matrix with unit columns $v_1, \dots, v_m$ satisfying
$A_{11}V = VD$ with  $D = \alpha_1 I_{m_1} \oplus \cdots \oplus \alpha_k I_{m_k}$ with $\alpha_1 > \dots > \alpha_k>0$ the distinct eigenvalues of $A_{11}$.
Let
$$\Omega_0 = \{\Gamma = \Gamma_1 \oplus \cdots \oplus \Gamma_k \in M_{m_1} \oplus \cdots \oplus M_{m_k}:
\Gamma \hbox{ is positive semidefinite} \},$$
$$\Omega_1 = \{\Gamma \in M_m:
D^{1/2}V^*(A_{11}A_{11}^* + A_{12}A_{12}^*)^{-1}VD^{1/2}\ge \Gamma \ge 0\},$$
and
$$\Omega_2 = \{\Gamma \in M_m: \Gamma  \ge V^*V\}.$$
The following proposition can be readily verified. Here we use the notation $X^{+}$ for the positive
semidefinite part of a Hermitian matrix $X$, i.e., $X^+ = (X +  \sqrt{X^2})/2$.

\begin{proposition} Let $G=[G_{ij}]$ be a Hermitian matrix, where $G_{ii}\in M_{m_i}$.
\begin{enumerate}
\item The projection of $G$ onto $\Omega_{0}$ is $G_{11}^{+}\oplus \cdots \oplus G_{kk}^{+}$.
\item The projection of $G$ onto $\Omega_{1}$ is $M-(M-G)^{+}$, where $M=D^{1/2}V^*(A_{11}A_{11}^* + A_{12}A_{12}^*)^{-1}VD^{1/2}$.
\item The projection of $G$ onto $\Omega_{2}$ is $(G-V^*V)^{+}+V^*V$.
\end{enumerate}
\end{proposition}
In the following algorithm, we create a sequence
\[\Gamma_0\longrightarrow \hat{\Gamma}_1\longrightarrow \Gamma_1 \longrightarrow \hat{\Gamma}_2 \longrightarrow \Gamma_2\longrightarrow \cdots\]
where $\Gamma_k\in \Omega_0$, $\hat{\Gamma}_{2k-1}\in \Omega_1$  and $\hat{\Gamma}_{2k}\in \Omega_2$ for all $k\geq 1$. This sequence converges to a solution $\Gamma \in \Omega_0\cap \Omega_1\cap \Omega_2$, provided $\Omega_0\cap \Omega_1\cap \Omega_2\neq \emptyset$; see \cite{AltProj}.

\begin{algorithm} For checking the existence of $\Gamma \in \Omega_0\cap \Omega_1\cap \Omega_2$.

\medskip
\noindent
Step 0. Set $k = 0$. Let $X=D^{1/2}V^*(A_{11}A_{11}^* + A_{12}A_{12}^*)^{-1}VD^{1/2}$
and $Y=V^*V$.

Partition  $X$  into $[X_{ij}]$ and $Y$ into $[Y_{ij}]$, both conformed to $D$.

Set $\Gamma_0=\frac{1}{2}\Big((X_{11}+Y_{11})\oplus\cdots \oplus (X_{kk}+Y_{kk})\Big)$. Go to Step 1.

\medskip\noindent
Step 1. Change $k$ to $k+1$, and set
\[\hat{\Gamma}_{k}=  \left\{\begin{array}{ll}
 X-(X-\Gamma_{k-1})^{+} & \mbox{ if $k$ is odd},\\
 (\Gamma_{k-1}-Y)^{+}+Y & \mbox{ if $k$ is even},
\end{array}\right. \]

 where $M_{+}$ denotes the positive part of $M$.

Partition $\hat{\Gamma}_{k}$ into $[G_{ij}]$ conformed to $D$
and let $\Gamma_{k}=G_{11}^{+}\oplus \cdots \oplus G_{kk}^{+}$.

If $error=\max(0,-\lambda_{\min}(\Gamma_{k}-Y))+\max(0,-\lambda_{\min}(X-\Gamma_k))\approx 0$, stop.

Otherwise, go to step $1$.

\end{algorithm}
Once we have $\Gamma$, we can set $U = V\Gamma^{-1/2}$, and
construct the two projections as shown in Theorem \ref{thm1}.
In particular, we can set $A = (I_p \oplus P)(I_p \oplus Q)$ with
\begin{equation}\label{pq}
P = {\small \begin{pmatrix} UU^*& 0 \cr 0 & 0_{n-p-m} \cr\end{pmatrix}}
\quad \hbox{ and } \quad
Q = {\small \begin{pmatrix} (U^*)^{-1}DU^{-1}& (UU^*)^{-1}A_{12} \cr
A_{12}^*(UU^*)^{-1} &  A_{12}^*(UDU^*)^{-1}A_{12}\cr\end{pmatrix}}.
\end{equation}

\medskip
We illustrate our Matlab program (available at http://cklixx.wm.edu/mathlib/Twoposcon.txt)
for checking whether a given matrix $A \in M_n$
is the product of two positive contractions in the following. Note that all numerical experiments were performed using Matlab 2015a on a Intel(R) Core(TM) i7-5500U CPU \@2.4GHz  with 8GB RAM and a 64-bit OS.

\begin{example}
Suppose $A={\small \begin{bmatrix}
A_{11} & A_{12} \\
0_{5} & 0_5
\end{bmatrix}}$, where \small
\[A_{11}={\small \begin{bmatrix}
0.125 & 0.0126  & 0.0033 & 0.024 & -0.0006\\
0 & 0.0625 & 0 & 0.012  & 0.0152\\
0 & 0 & 0.0625 & 0.0025 & 0.0453\\ 0 & 0 & 0 & 0.2 & 0\\
0 & 0  & 0  & 0  & 0.2
\end{bmatrix}} \mbox{ and } A_{12}={\small \begin{bmatrix}
0.0658 & 0.0218  & 0.0031 & 0.05  & -0.0033 \\ 0.0218 & 0.113  & -0.0107 & -0.0120 &  0.0098\\
0.0031& -0.0107 & 0.0418 & 0.0048 & -0.0409\\
0.0500 & -0.012  & 0.0048 & 0.1103 & 0.0037\\
-0.0033 & 0.0098  & -0.0409 & 0.0037 & 0.128
\end{bmatrix}}.\]
We set
\[V\approx \begin{bmatrix}
1 &	-0.1976 & -0.0507 &	-0.3169 & 	-0.0169 \\
0 & 0.9803 & -0.0102 & -0.0824 & 	-0.1026\\
0 & 0 & 0.9987 & -0.0172 & -0.3108 \\
0 & 0 & 0 & -0.9447 & 0.0203\\
0 & 0 & 0 & 0 & -0.9445
\end{bmatrix},\]
which has unit columns and satisfies $A_{11}V=V{\rm diag}(0.125,0.0625,0.0625,0.2,0.2)$;
the second and third columns of $V$ are orthogonal and the fourth and fifth columns
are orthogonal.

\medskip\noindent
\rm Using our Matlab program, we obtain
$U=V\Gamma^{-\frac{1}{2}}$, where
\[\Gamma=\begin{bmatrix}
3.4737 & 0 & 0 & 0 & 0\\
0 & 2.3344 & 0.0216  & 0 & 0\\
0 & 0.0216 & 2.9472 & 0 & 0\\
0 & 0 & 0 & 2.1257 & -0.2132\\
0 & 0 & 0 &-0.2132 & 1.6425
\end{bmatrix}. \]
Defining $P$ and $Q$ as in equation (\ref{pq}), we get that $\lambda_1(P)=s_1^2(U)=0.7024$ and $\lambda_1(Q)=1$. Note that $\Gamma$ is obtained using alternating projection method after 79 iterations done in approximately 0.085 seconds with $errror=||PQ-A||=4.3774\times 10^{-14}$.
\normalsize
\end{example}

\begin{example} \small
Suppose
\[A_{11}=\begin{bmatrix}
0.1 &  0.0244 & 0.026 & 0.0167 & 0.0114 & 0.0014 & 0.0674\\
0 & 0.2 & 0.0176 & 0.0251 & 0.0345 & 0.0122 & 0.0088\\
0 & 0 & 0.3 & 0 & 0.0072 & 0.0119 & 0.0166\\
0 & 0 & 0 & 0.3 & 0.0093 & 0.0007 & 0.0099\\
0 & 0 & 0 & 0 & 0.4 & 0 & 0\\
0 & 0 & 0 & 0 & 0 & 0.4 & 0\\
0 & 0 & 0 & 0 & 0 & 0 & 0.4
\end{bmatrix}\]
and
\[A_{12}=\begin{bmatrix}
0.098 & 0.0157 & -0.0315 & 0.0033 & -0.04 & -0.0196 & 0.0171\\
0.0157 & 0.0545 & -0.0366 & 0.0302 & 0.0081 & 0.0003 & 0.004\\
-0.0315 & -0.0366 & 0.1246 & -0.0449 & -0.0005 &
0.0232 & -0.0047\\
0.0033 & 0.0302 & -0.0449 & 0.1025 & -0.0193 & -0.031 & 0.0191\\
-0.04 & 0.0081 & -0.0005 & -0.0193 & 0.1285 & 0.0038 & -0.0504\\
-0.0196 & 0.0003 & 0.0232 & -0.031 & 0.0038 & 0.07790 & -0.0192\\
0.0171 & 0.004 & -0.0047 & 0.0191 & -0.0504 & -0.0192 & 0.0895
\end{bmatrix}. \]
We let \[V=\begin{bmatrix}
1 & -0.2373 & -0.1475 & -0.1015 & -0.0632 & -0.0196 & -0.2348\\
0 & -0.9714 & -0.1713 & -0.2329 & -0.1858 & -0.0673 & -0.0569\\
0 & 0 & -0.9741 & 0.0563 & 	-0.0702 & -0.1162 & 	-0.1512 \\
0 & 0 & 0 & -0.9656 & -0.0910 & -0.0052 & 	-0.0896 \\
0 & 0 & 0 & 0 & -0.9738 & 0.023 & 0.0454 & \\
0 & 0 & 0 & 0 & 0 & -0.9905 & 0.0278\\
0 & 0 & 0 & 0 & 0 & 0 & -0.9528
\end{bmatrix}.\]
\rm
Using our Matlab program, we obtain
\[\Gamma= [2.9099] \oplus [2.592]\oplus \begin{bmatrix}
1.9048 &	0.1063\\
0.1063 & 1.866
\end{bmatrix}\oplus \begin{bmatrix}
1.6447 & 0.0046 & 0.0768 \\
0.0046 & 1.6923 &	0.0215\\
0.0768 & 0.0215 & 1.5846
\end{bmatrix} \]
after 59 iterations (approximately 0.075 seconds) with a $1.227\times 10^{-16}$ error. The positive semidefinite matrices $P$ and $Q$ defined in equation (\ref{pq}) will have largest eigenvalues $0.8309$ and $1$, respectively.
\normalsize
\end{example}

\begin{example} Let

\medskip
\centerline{
$A=\begin{bmatrix}
A_{11} & 0 \\ 0 & 0
\end{bmatrix}$ and $B=\begin{bmatrix}
B_{11} & 0 \\ 0 & 0
\end{bmatrix}$, where
$A_{11}=\begin{bmatrix}
0.5 & 0.09429 \\
0 & 0.3
\end{bmatrix}$ and
$\quad B_{11}=\begin{bmatrix}
0.5 & 0.0943 \\
0 & 0.3
\end{bmatrix}$.}

\medskip\noindent\rm
 It follows from \cite{LT} that $A$ is a product of two contractions and $B$ is not. Notice that $A$ and $B$ are very close to each other.

For $A$, we ran the alternating projection algorithm and  obtained $\Gamma=\mbox{diag}(1.2759,1.6591)$ after 66321 iterations (48.26 seconds). We also get $||PQ-A||\approx 1.4778\times 10^{-16}$ and $\lambda_1(Q),\lambda_1(P)\approx 1$. \\
Meanwhile, for $B$, after running 100,000 iterations (69.06 seconds) of the algorithm, we see that
the  values $\max(0,-\min({\rm eig}(M-\Gamma)))$ and $\max(0,-\min({\rm eig}(\Gamma-V^*V)))$
starts to alternate back and forth from $8.5\times 10^{-5}$ to $8.52925\times 10^{-5}$.
\end{example}

\newpage
\bigskip\noindent
{\large \bf Acknowledgment}

The authors would like to thank Professor T. Ando for some helpful discussion.
We also thank the referee for his/her comments and suggestions that help improve the paper. 
Li is an
honorary professor of the University of Hong Kong and the Shanghai University. His research was
supported by USA NSF grant DMS 1331021, Simons Foundation Grant 351047, and NNSF of China
Grant 11571220.
The research
of Wang was supported by the Ministry of Science and Technology of the Republic of China under
project MOST 104-2115-M-009-001

\end{document}